\documentclass[sn-mathphys]{sn-jnl}
\jyear{2021}
\usepackage{microtype}
\theoremstyle{thmstyleone}%
\newtheorem{theorem}{Theorem}
\newtheorem{proposition}[theorem]{Proposition}%

\newtheorem{lemma}[theorem]{Lemma}

\newtheorem{counterexample}[theorem]{Counterexample}

\numberwithin{equation}{section}

\theoremstyle{thmstyletwo}%
\newtheorem{remark}{Remark}%

\theoremstyle{thmstylethree}%
\newtheorem{definition}{Definition}%

\newcommand{\mC}{\ensuremath{\mathbb{C}}}
\newcommand{\mD}{\ensuremath{\mathbb{D}}}

\newcommand{\mN}{\ensuremath{\mathbb{N}}}
\newcommand{\mR}{\ensuremath{\mathbb{R}}}

\begin{document}

\title{Counterexamples to Bloch's Principle and its Converse in Several Complex Variables}


\author[1]{\fnm{Kuldeep} \sur{Singh Charak}}\email{kscharak7@rediffmail.com}

\author*[2]{\fnm{Rahul} \sur{Kumar}}\email{rktkp5@gmail.com}


\affil[1]{\orgdiv{Department of Mathematics}, \orgname{University of Jammu}, \orgaddress{\city{Jammu}, \postcode{180006}, \state{Jammu and Kashmir}, \country{India}}}

\affil*[2]{\orgdiv{Department of Mathematics}, \orgname{University of Jammu}, \orgaddress{\city{Jammu}, \postcode{180006}, \state{Jammu and Kashmir}, \country{India}}}


\abstract{In this paper, besides a counterexample to Bloch's principle, normality criteria leading to counterexamples to the converse of Bloch's principle in several complex variables are proved. Some Picard-type theorems and their corresponding normality criteria in $\mC^n$ are also obtained.}

\keywords{Normal families, Picard type theorems, Bloch's principle}


\pacs[MSC Classification]{32A19; 32A10}

\renewcommand{\thefootnote}{\fnsymbol{footnote}}
\footnotetext{The work of the second author is supported by CSIR, India (File No. 09/100(0217)/2018-EMR-I).}

\maketitle



This paper is divided into the following four sections:
\begin{enumerate}
    \item Introduction and Auxiliary Results;
    \item Counterexamples to Bloch's Principle and its Converse;
    \item Picard-type Theorems and their Corresponding Normality Criteria;
    \item Extension of some known Results from $\mC$ to $\mC^n$.
\end{enumerate}

\section{\textbf{Introduction and Auxiliary Results }}

Let $D$ to be a domain in $\mC^n$ and $\mathcal F$ to be a family of holomorphic functions in $D$.  $\mathcal F$ is said to be normal in $D$ if each sequence in $\mathcal F$ contains a subsequence that converges locally uniformly in $D.$ The aim of this paper is to obtain normality criteria leading to counterexamples to the converse of Bloch's principle in several complex variables. Bloch's principle in one complex variable is extensively studied and the reader may refer to \cite{Robinson, Rubel, Zalcman, Bergweiler, Li, Lahiri, Charak1, Charak2, Charak3}, in order. Bloch's principle states that a family  $\mathcal F$ of holomorphic functions of several complex variables in a domain $D\subset \mC^n$ satisfying a certain property $\bf P$ in $D$ is likely to be normal in $D$ if any entire function  possessing the property $\bf P$ in $\mC^n$ reduces to a constant. In one complex variable neither Bloch's principle  nor its converse holds in general, for example one may refer to \cite{Rubel, Lahiri}. For the construction of counterexamples to the converse, one needs to find a normality criterion with a certain property and then find out a non-constant entire function satisfying this property in $\mC^n$; this sounds interesting particularly from the experience of one complex variable case and this is what we have explored in the present paper in several complex variables case too.\\
 We shall denote by $\mathcal {H}(D)$, the class of holomorphic functions $f:D\rightarrow \mC$,  where $D\subset \mC^n$ is a domain, and the open unit ball in $\mC^n$ shall be denoted by $\mathbb{D}^n :=\{z\in\mC^n:\|z\|<1\}.$\\
For every $\psi\in \mathcal{C}^2(D),$ at each point $z\in D$ we define a Hermitian form
\begin{equation}\label{eqA}
L_z(\psi,~v):= \sum\limits_{k,l=1}^{n} \frac{\partial^2 (\psi)}{\partial z_k\partial \bar{z_l}}(z)v_k\bar{v_l}
\end{equation} 
and is called the {\it Levi form } of the function $\psi$ at $z.$\\
The spherical derivative of $f\in \mathcal{H}(D)$ (see \cite{Dov1})  is defined as 
\begin{equation}\label{eqB}
f^{\#}(z):= \sup\limits_{\|v\|=1}\sqrt{L_z(\log(1+\lvert f\rvert^2),~v)}.
\end{equation} 
Since $L_z(\log(1+\lvert f\rvert^2),~v)\geq 0,  \  f^{\#}(z)$ given by (\ref{eqB}) is well defined  and for $n=1$, $\eqref{eqB}$ reduces to 
$$f^{\#}(z):= \frac{\lvert f^\prime(z)\rvert}{1+\lvert f(z)\rvert^2}$$
which is the spherical derivative of a holomorphic function of one complex variable. Thus (\ref{eqB}) gives the natural extension of the spherical derivative to $\mC^n.$\\
Further, from  (\ref{eqB}), one easily finds that
\begin{equation}\label{eqC}
f^{\sharp}(z) = \frac{\sup\limits_{\|v\|=1} \lvert (Df(z), v)\rvert}{1+\lvert f(z)\rvert^2},
\end{equation}
where
$$(Df(z),v)=\sum _{j=1}^{n}\frac{\partial f(z)}{\partial z_j}.v_j$$

Following known results shall be required to accomplish the proofs of main results of this paper, spread over various sections:\\

\medskip

\noindent {\bf Marty's Theorem in $\mC^n$ (\cite{Dov1})}: {\it A family $\mathcal{F}\subset \mathcal{H}(D)$ is normal in $D$ if and only if for each compact set $K\subset D $ there exists a constant $M=M(K)$ such that  $f^{\#}(z)\leq M $ for all  $f\in\mathcal{F}, \ z\in K.$}

\medskip

\noindent {\bf Zalcman's Lemma in $\mC^n$} (\cite{Dov2}): 
{\it A family $\mathcal{F}\subseteq \mathcal{H}(D)$ is not normal  at a point $z_0\in D$ if and  only if for each $\alpha\in (-1,\infty)$ there exist sequences $\{f_j\}\subset \mathcal{F}, ~\{z_j\}\subset D : z_j\to z_0, \mbox{ and } ~\{r_j\}\subset (0,1]:r_j\to 0$ such that the zoomed sequence
$$g_j(z):= r_j^{\alpha}f_j(z_j+r_jz)$$ 
converges locally uniformly to a non-constant entire function $g$ in $\mC^n$ satisfying $g^{\sharp}(z)\leq g^{\sharp}(0)=1.$}

\begin{definition}
Let $f\in \mathcal{H}(\mC^n)$   and $ z= (z_1,z_2,\ldots,z_n)\in\mC^n.$ The total derivative of $f$, which we shall denote by $D_0f,$ is defined as 
$$D_0f(z)=\displaystyle\sum_{j=1}^{n}z_jf_{z_j},$$ where $f_{z_j}$ is the partial derivative of $f$ with respect to $z_j.$
\end{definition}

Liu and Cao \cite{Liu} obtained an extension of Zalcman's lemma concerning the total derivative in several complex variables as follows:
\begin{lemma}\label{lemma1}
Let $\mathcal{F}\subset \mathcal{H}(\mathbb{D}^n)$  and suppose that $f(z)\neq 0$ for every $f\in\mathcal{F}$ and for all $z\in\Delta_\delta(0)$ for some $\delta \in (0,1)$ and  that $f(z)= 0 \implies \lvert D_0f(z)\rvert\leq A$ for some $A>0.$ If $\mathcal{F}$ is not normal in $\mathbb{D}^n,$ then for all $k\in (-1,  1],$ there exist a real number $r\in (0,  1],$ and sequences  $\{z_j\}\subseteq \mathbb{D}^n :0<\|z_j\|<r, \  \{f_j\}\subseteq\mathcal{F},$ and $\{\rho_j\}\subset (0,1]: \rho_j \to 0$ such that 
$$g_j(\zeta)= \rho_j^{-k}f_j(z_je^{\rho_j\zeta}), ~(\zeta\in\mC)$$ 
converges locally uniformly to a non-constant entire function $g$ in $\mC$  satisfying  $g^\sharp(\zeta)\leq g^\sharp(0)= A+1.$
\end{lemma}
\begin{definition}
For $f\in\mathcal{H}(\mC^n) ~\mbox{and}~ a\in\mC^n,$ we can write $f(z) = \sum\limits_{i=0}^{\infty}P_i(z-a),$ where $P_i(z)$ is either a homogeneous polynomial of degree $i$ or identically zero. The zero multiplicity of $f$ at $a$ is defined as $V_f(a) = min\{i:P_i\neq 0\}.$
\end{definition}
We have a kind of variant of Lemma \ref{lemma1}:
\begin{lemma}\label{thm ZP}
Let $\mathcal{F}\subset \mathcal{H}(\mD^n).$  For $ f\in\mathcal{F},$ suppose that for each $a\in f^{-1}(\{b\})\cap \mD^n,$ zero multiplicity of $f ~\mbox{at}~ a$ is at least $k.$ Then if $\mathcal{F}$ is not normal, then for each $\alpha \in (0,  k),$ there exist a real number $r\in (0,  1)$ and sequences  $\{z_j\}\subset \mD^n: \|z_j\|<r, \ \{f_j\}\subset \mathcal{F},$ and $\{\rho_j\}\subset (0,1]:\rho_j\to 0$ 
such that 
$$g_j(\zeta)= \rho_j^{-\alpha}f_j(z_j+\rho_j\zeta)$$ 
converges locally uniformly to a non-constant entire function $g$ in $\mC^n$ satisfying  $g^{\sharp}(\zeta)\leq g^{\sharp}(0) = 1.$
\end{lemma}
The proof of Lemma \ref{thm ZP} is exactly on the lines of proof of Theorem $5$ in \cite{Charak4} with Lemma $1$ there replaced by  
\begin{lemma}\label{l1}
 Let $f\in \mathcal{H}(\mD^n),$  for each $a\in f^{-1}(\{b\})\cap \mD^n,$ zero multiplicity of $f ~\mbox{at}~ a$ be at least $k, $  and $\alpha \in (-1, k).$ Let $\Omega :=\{z\in \mC^n :\|z\|<r<1\}\times (0,  1]$ and define a function $F:\Omega \rightarrow \mR$ as   

$$F(z,t)=\frac{\left(1-\frac{\|z\|^2}{r^2}\right)^{1+\alpha}t^{1+\alpha}f^{\#}(z)\left(1+\lvert f(z)\rvert^2\right)}{\left(1-\frac{\|z\|^2}{r^2}\right)^{2\alpha}t^{2\alpha}+\lvert f(z)\rvert^2}.$$
If $F(z,1)>1$ for some $z:\|z\|<r$, then there exist $z_0\in \mC^n :\|z_0\|<r<1$ and $t\in (0,  1)$ such that

$$\sup\limits_{\|z\|<r} F(z,t_0)=F(z_0,t_0)=1.$$
\end{lemma}
Since the proof of Lemma \ref{l1} is a modification of the proof of Lemma $1$ in \cite{Charak4}, here we  only need to establish that  
\begin{equation}\label{eq1}
\lim\limits_{\left(1-\frac{\|z\|^2}{r^2}\right)t\to 0}F(z,t) = 0,
\end{equation}
as a significant amount of modifications is required to establish this limit.

\smallskip

We may assume that $z_j\to z_0,$ where $\|z_0\|\leq r.$\\
If $f(z_0)\neq 0,$ then

$$\lim\limits_{j\to\infty} F(z_j,t_j) \leq 0.$$ 
Suppose $f(z_0)= 0.$ Then $f(z)= \sum_{m=s}^\infty P_m(z-z_0), ~s\geq k,$ where $P_m(z)$ is either a homogeneous polynomial of degree $m$ or identically zero.

\begin{eqnarray*}
\lvert f(z_j)\rvert &=& \lvert p_s(z_j-z_0)+p_{s+1}(z_j-z_0)+\ldots \rvert\\
	       &=& \|z_j-z_0\|^s\left[\frac{\lvert p_s(z_j-z_0)+p_{s+1}(z_j-z_0)+\ldots\rvert}{\|z_j-z_0\|^s}\right]\\
				 &=& \lvert a_s \rvert \|z_j-z_0\|^s[1+o(1)].
\end{eqnarray*}

Similarly, we find that
$$\left\lvert \frac{\partial}{\partial z_i}f(z_j)\right\rvert= \lvert b_{s_i}\rvert \|z_j-z_0\|^{s-1}[1+o(1)].$$
We may assume that 
$$\lim\limits_{j\to\infty}\frac{\left(1-\frac{\|z_j\|^2}{r^2}\right)t_j}{\|z_j-z_0\|} = B$$ 
exists. If $B\geq 1,$ then
\begin{eqnarray*}
\lim\limits_{j\to\infty}F(z_j,t_j)  &\leq& \lim\limits_{j\to\infty}\left(1-\frac{\|z_j\|^2}{r^2}\right)^{1-\alpha}t_j^{1-\alpha}f^{\#}(z_j)\left(1+\lvert f(z_j)\rvert^2\right)\\
	                                  &\leq& \frac{\|\left(\frac{\partial}{\partial z_1}f(z_j),\ldots,\frac{\partial}{\partial z_n}f(z_j)\right)\|}{\left(1-\frac{\|z_j\|^2}{r^2}\right)^{\alpha-1}t_j^{\alpha-1}}\\
																		&\leq&  \lim\limits_{j\to\infty} \frac{\lvert \frac{\partial}{\partial z_1}f(z_j)\rvert+\ldots+\lvert \frac{\partial}{\partial z_n}f(z_j)\rvert}{\left(1-\frac{\|z_j\|^2}{r^2}\right)^{\alpha-1}t_j^{\alpha-1}}\\
																		&\leq& \lim\limits_{j\to\infty} \frac{\lvert b_{s_1}\rvert\|z_j-z_0\|^{s-1}+\ldots \lvert b_{s_n}\rvert\|z_j-z_0\|^{s-1}}{\left(1-\frac{\|z_j\|^2}{r^2}\right)^{\alpha-1}t_j^{\alpha-1}}\\
																		&\leq& \lim\limits_{j\to\infty} \frac{n \max_{i=1}^{n} \lvert b_{s_i}\rvert\|z_j-z_0\|^{s-1}}{\left(1-\frac{\|z_j\|^2}{r^2}\right)^{\alpha-1}t_j^{\alpha-1}}\\
																		&=& \lim\limits_{j\to\infty} \frac{n \max_{i=1}^{n} \lvert b_{s_i} \rvert\|z_j-z_0\|^{s-\alpha}\|z_j-z_0\|^{\alpha-1}}{\left(1-\frac{\|z_j\|^2}{r^2}\right)^{\alpha-1}t_j^{\alpha-1}}\\
																		&=& 0
\end{eqnarray*}
	If $B\leq 1,$ then
	$$F(z_j,t_j)= \frac{\left(1-\frac{\|z_j\|^2}{r^2}\right)^{\alpha}t_j^{\alpha}\lvert f(z_j)\rvert}{\left(1-\frac{\|z_j\|^2}{r^2}\right)^{2\alpha}t_j^{2\alpha}+\lvert f(z_j)\rvert^2}.\left(1-\frac{\|z_j\|^2}{r^2}\right)t_j\frac{f^{\#}(z_j)\left(1+\lvert f(z_j)\rvert^2\right)}{\lvert f(z_j)\rvert}.$$

 \medskip
 
The first factor on the right hand is obviously bounded by $\frac{1}{2}.$\\
Moreover, if $B\neq 0,$ then we have
	\begin{eqnarray*}
	\lim\limits_{j\to\infty}\frac{\left(1-\frac{\|z_j\|^2}{r^2}\right)^\alpha t_j^\alpha \lvert f(z_j)\rvert}{\left(1-\frac{\|z_j\|^2}{r^2}\right)^{2\alpha}t_j^{2\alpha}+\lvert f(z_j)\rvert^2} &\leq&
\lim\limits_{j\to\infty}\frac{B^\alpha \|z_j-z_0\|^\alpha\lvert a_s\rvert\|z_j-z_0\|^s}{B^{2\alpha}\|z_j-z_0\|^{2\alpha}}\\
	&=& \lim\limits_{j\to\infty}\frac{\lvert a_s\rvert\|z_j-z_0\|^{s-\alpha}}{B^\alpha} \\
	&=& 0
	\end{eqnarray*}
	
On the other hand,
	\begin{eqnarray*}
	\lim\limits_{j\to\infty}\left(1-\frac{\|z_j\|^2}{r^2}\right)t_j\frac{f^{\#}(z_j)\left(1+\lvert f(z_j)\rvert^2\right)}{\lvert f(z_j)\rvert} &\leq& \left(1-\frac{\|z_j\|^2}{r^2}\right)t_j.\frac{n\max_{i=1}^n\lvert b_{s_i}\rvert\|z_j-z_0\|^{s-1}}{\lvert a_s\rvert\|z_j-z_0\|^s}\\
	&=& \lim\limits_{j\to\infty}\left(1-\frac{\|z_j\|^2}{r^2}\right)t_j.\frac{n\max_{i=1}^n\lvert b_{s_i}\rvert}{\lvert a_s\rvert\|z_j-z_0\|}\\
	&=& B\frac{n\max_{i=1}^n\lvert b_{s_i} \rvert}{\lvert a_s\rvert}
	\end{eqnarray*}
	It follows that $\lim\limits_{j\to\infty}F(z_j,t_j) = 0$.

\medskip

\bigskip

Let $f$ be meromorphic function in $\mC$ and $ a\in \mC_\infty.$ Then $a$ is called totally ramified value of $f$  if $f-a$ has no simple zeros.
Following result known as {\it Nevanlinna's Theorem } (see \cite{Bergweiler}) plays a crucial role in the  proofs of Theorem \ref{BP}, Theorem \ref{P1} and Theorem \ref{P2} to follow:
\begin{theorem}\label{thmN}
Let $f$ be a non-constant meromorphic function $a_1,...,a_q\in \mC_\infty \mbox{ and } m_1,...,m_q \in \mN.$ Suppose that all $a_j$-points of $f$ have multiplicity  at least $m_j,  \mbox{ for }~ j=1,...,q.$ Then 
$$\sum\limits_{j=1}^{q}\left(1-\frac{1}{m_j}\right) \leq 2.$$
\end{theorem}

\section{Counterexamples to Bloch's Principle and its Converse}
\subsection{A Counterexample to Bloch's Principle}
If $D\subset \mC^n$ and $f\in \mathcal{H}(D)$ satisfies a certain property $\bf P$ in $D, $ then we write it as $(f,D) \in \bf P.$ Recall that Bloch's principle states that a family $\mathcal{F}\subset \mathcal{H}(D): (f,D)\in \bf P$ is likely to be normal in $D$ if $(f, \mathbb{C}^n)\in \bf P$ implies that $f$ is a constant function. As in the case of one complex variable, Bloch's principle fails to hold in $\mathbb{C}^n$ as well.\\
Let $\bf P$ be the property of holomorphic function $f$ in $D$ defined as  
$$P(f)(z)= \displaystyle\prod_{i=1}^{n}\left(\frac{\partial f(z)}{\partial z_i}-a\right)\left(\frac{\partial f(z)}{\partial z_i}-b\right)\left(\frac{\partial f(z)}{\partial z_i}-f(z)\right)$$ omits the value zero. \\
\begin{proposition}
If $f\in\mathcal{H}(\mC^n) : (f,\mC^n)\in \bf P,$ then $f\equiv$constant.
\end{proposition}

\textbf{Proof:}   Suppose $f\in\mathcal{H}(\mC^n)$ possesses the property $\bf P$. Let 
$$w= (a_1,\ldots, a_n), ~w^{\prime}= (b_1,\ldots,b_n) \in\mC^n$$ 
and define 
$$h_i(z_i)= f(b_1,\ldots,b_{i-1},z_i,a_{i+1},\ldots,a_n), ~i= 1,2,\ldots n.$$
Then $\frac{d}{d z_i}h_i(z_i)$ omits the values $a$ and $b$ and so $\frac{d}{d z_i}h_i(z_i)= c_i, $ a constant. This implies that $h_i(z_i)= c_i z_i+d_i, $ but  $\frac{d}{d z_i}h_i(z_i)-h_i(z_i)$ omits zero and so $c_i-c_iz_i-d_i \neq 0$ which implies $c_i = 0.$ This shows that each $h_i(z_i)$ is constant and hence $f$ is constant.$\Box$

\begin{counterexample} Following counterexample to the Bloch's principle is constructed from the ideas of Rubel\cite{Rubel}.\\
    
Take $r$ sufficiently small so that 
$$\{(z_1,\ldots,z_n):z_1+\ldots +z_n=1\}\cap D^n(0,r)= \emptyset .$$
Consider 
$$\mathcal{F}= \{f_m(z_1,\ldots,z_n)= m(z_1+\ldots +z_n) :(z_1,\ldots,z_n)\in D^n(0,r), m\geq 3\}.$$
Then 
$$\frac{\partial f(z)}{\partial z_i} = m\neq 1, 2.$$  
Also, $$\frac{\partial f(z)}{\partial z_i}-f(z)= m(1-(z_1+\ldots+z_n))\neq 0 ~~\forall z=(z_1,\ldots,z_n)\in D^n(0,r).$$ 
Thus each $f_m$ satisfies the property $\bf P$ in $D^n(0,r)$ with $a=1 ~\mbox{and}~ b=2$ but $\mathcal{F}$ is not normal in $D^n(0,r)$ showing that Bloch's principle does not hold in general. 
\end{counterexample}

\subsection{Counterexamples to the Converse of Bloch's Principle}

 The converse of Bloch's principle precisely states that if 
 $$\mathcal{F}:=\{f\in \mathcal{H}(D): (f,D) \in {\bf P}\}$$
 is normal in $D,$ then  any $f\in\mathcal{H}(\mC^n): (f, \mC^n)\in {\bf P}$ reduces to a constant.\\
 Let the property {\bf P} satisfied by a function $f\in \mathcal{H}(D)$ be defined as follows:
 \begin{itemize}
\item [(i)] there exist three distinct complex numbers $a_1, a_2, a_3$ and three complex numbers $b_1, b_2, b_3$ satisfying $f(z)= a_i \implies D_0f(z)= b_i;$
\item [(ii)] $f(z)\neq0 ~\forall ~z\in \Delta_\delta(0) ~\mbox{for some}~ \delta>0.$
\end{itemize}
This property leads to the following normality criterion:
\begin{theorem}\label{BP}
 A subfamily $\mathcal{F}$ of $\mathcal{H}(D)$ is normal in $D$ if $(f, D)\in {\bf P}$ for all $f\in \mathcal{F}.$
\end{theorem}
\textbf{Proof:}
Since normality is a local property, we may assume $D$ to be the open unit ball $\mathbb{D}^n:= \{z\in\mC^n:\|z\|<1\}.$ Suppose $\mathcal{F}$ is not normal in $\mathbb{D}^n.$ Then by Lemma \ref{lemma1}, there exist $r\in (0,  1), ~z_j\in \mD^n :0<\|z_j\|<r, ~\{f_j\}\subseteq \mathcal{F}, ~\{\rho_j\}\subset (0,  1]:\rho _j\to 0$ such that $$g_j(\zeta)= f_j(z_je^{\rho_j\zeta})$$ 
converges locally uniformly to non-constant entire function $g$ in $\mC.$ Let $g(\zeta_0)= a_i.$ Then by Hurwitz's theorem there exists a sequence $\{\zeta_j\}$ converging to $\zeta_0$ such that for sufficiently large $j, ~g_j(\zeta_j)= a_i$ which implies $f_j(z_je^{\rho_j\zeta_j})= a_i$ and so by given hypothesis $D_0f_j(z_je^{\rho_j\zeta_j})= b_i.$ Therefore,

$$g^\prime(\zeta_0)= \lim_{j\to\infty}g_j^\prime(\zeta_j)= \lim_{j\to\infty}\rho_j D_0f_j(z_je^{\rho_j\zeta_j})= 0.$$

Thus each $a_i$-points of $g$ has multiplicity at least two which contradicts Theorem \ref{thmN}. $\Box$

\medskip

Theorem \ref{BP} leads to a counterexample to the converse of Bloch's principle.
\begin{counterexample}
Consider the function
$$f(z_1,z_2,\ldots,z_n) = z_1+z_2+\ldots+z_n+4,$$ 
$a_1= 1, \  a_2= 2, \  a_3= 3$ and $b_1=-3, \  b_2=-2, \  b_3= -1.$\\
 Then $f(z)= a_j \implies D_0f(z)= b_j ~(j= 1,2,3).$ Also $f(z)\neq 0 ~\forall ~z\in \Delta_\delta(0) ~\mbox{for some}~ \delta>0.$
Thus there exists a non-constant entire function in  $\mC^n$ satisfying the property $\bf{P}$ and by Theorem \ref{BP},
$\mathcal{F}:= \{f \in \mathcal{H}(D):(f,D)\in \bf{P}\}$ is normal in $D$. 
\end{counterexample}

The following weak version of Lappan's normality criterion in several complex variables(see \cite{Charak4},Theorem 1)  also leads to a counterexample to the converse of Bloch's principle:
\begin{theorem}\label{thm001} Let    $D\subseteq \mathbb{C}^{n}$ and  $\mathcal{F}\subset \mathcal{H}(D).$  Let $E$ be a subset of  $\mC$ containing at least three points such that for each compact set $K\subset D,$ there exists a positive constant $M=M(K)$ for which 
\begin{equation}\label{eq001}
\sup\limits_{\left\|v\right\|=1}\lvert (Df(z), v)\rvert \leq M \mbox{ whenever }  f(z)\in E, ~ z\in K, ~ f\in\mathcal{F}.
\end{equation}
Then  $\mathcal{F}$ is normal in $D.$
\end{theorem}
The proof of Theorem \ref{thm001} is obtained on the lines of the proof of Theorem $1$ in \cite{Charak4} hence omitted.

\begin{counterexample}
 Let $D\subset \mC^n$ be a domain. For $f\in \mathcal{H}(D),$  $(f, D)\in {\bf P}$ if there exists a set $E\subset \mC$ containing at least three points such that for each compact subset $K\subset D,$ there exists a positive constant $M=M(K)$ for which 
\begin{equation}
\sup\limits_{\left\|v\right\|=1}\lvert (Df(z), v)\rvert\leq M \mbox{ whenever }  f(z)\in E, ~ z\in K.
\end{equation} 
Consider the function 
$$f(z_1,z_2,\ldots,z_n) = z_1+z_2+\ldots+z_n.$$ Then 
$$\sup\limits_{\left\|v\right\|=1}\lvert (Df(z), v)\rvert\leq \sup\limits_{\left\|v\right\|=1}[\lvert v_1\rvert+\lvert v_2\rvert+\cdots+\lvert v_n\rvert]\leq n.$$\\
Thus $f$ is a non-constant entire function in  $\mC^n$ satisfying the property ${\bf P}$ in $\mC^n$ and by Theorem \ref{thm001},
$\mathcal{F}:= \{f \in \mathcal{H}(D):(f,D)\in {\bf P}\}$ is normal in $D$. This is a another counterexample to the
converse of Bloch's principle in several complex variables.
\end{counterexample}

\section{Picard-type Theorems and their Corresponding Normality Criteria}
By Bloch's principle one can expect to have a normality criterion corresponding to every Picard-type theorem. Here we prove a couple of Picard-type theorems and obtain the corresponding normality criteria.\\
In the discussions to follow, we shall denote the operator $\left(\frac{\partial^k}{\partial z_1^k},\frac{\partial^k}{\partial z_2^k},\ldots,\frac{\partial^k}{\partial z_n^k}\right)$ by $D^k.$
\begin{theorem}\label{P0}
 Let $f\in \mathcal{H}(\mC^n)$ and $k$ be a positive integer. If $f(z) \neq 0, ~(D^kf(z),v) \neq 1 ~\mbox{for all}~  ~v=(v_1,...,v_n):\|v\|=1,$ then $f$ is constant.\end{theorem}
\textbf{Proof:}
Suppose $f$ is non-constant.  Let 
$$w= (a_1,\ldots, a_n), ~w^{\prime}= (b_1,\ldots,b_n) \in\mC^n$$
and define 
$$h_i(z_i):= f(b_1,\ldots,b_{i-1},z_i,a_{i+1},\ldots,a_n), ~i= 1,2,\ldots n.$$
By given hypothesis, it follows that $h_i(z)\neq 0 ~\mbox{for all}~ z\in\mC.$

\smallskip

\noindent 
\textbf{Claim:} $\frac{\partial^k}{\partial z_i^k}f(b_1,\ldots,b_{i-1},z_i,a_{i+1},\ldots,a_n) \not\equiv 1$

\smallskip

\noindent
Suppose $$\frac{\partial^k}{\partial z_i^k}f(b_1,\ldots,b_{i-1},z_i,a_{i+1},\ldots,a_n) \equiv 1.$$ 
Then $h_i(z)$ is a polynomial of degree $k$ and so $h_i(a) = 0$ for some $a\in\mC,$ a contradiction and this establishes the claim.\\
Since $h_i(z)\neq 0 ~\mbox{for all}~ z\in\mC.$ So, $\frac{\partial^{k}}{\partial z_i^k}h_i(a)= 1 ~\mbox{for some}~ a\in\mC$ which implies 
$$\frac{\partial^k}{\partial z_i^k}f(b_1,\ldots,b_{i-1},a,a_{i+1},\ldots,a_n) = 1.$$
which implies that
$$\sum\limits_{j=1}^{n} \frac{\partial^k }{\partial z_j^k}f(b_1,\ldots,b_{i-1},a,a_{i+1},\ldots,a_n)e_j = 1$$
where $e_j=
\left\{\begin{array}{cc}
1, & ~\mbox{for}~ j=i\\
0, & ~\mbox{for}~ j\neq i.\end{array}\right.$

\noindent 
That is,
$$(D^kf(b_1,\ldots,b_{i-1},a,a_{i+1},\ldots,a_n),v) = 1, ~\mbox{where}~ v = (0,\ldots,1,\ldots, 0),$$ which is a contradiction.    $\Box$

\medskip

The corresponding normality criterion of Theorem \ref{P0} is
\begin{theorem}\label{thm M}
Let $D$ be a domain in $\mC^n$. Let $\mathcal{F}$ be a family of holomorphic functions in a domain $D$ and $k$ be a positive integer. If $f(z) \neq 0, ~(D^kf(z),v) \neq 1 $ for all $f\in\mathcal{F}, ~v=(v_1,...,v_n):\|v\|=1,$ then $\mathcal{F}$ is normal in $D.$
\end{theorem}
\textbf{Proof:}
Suppose $\mathcal{F}$ is not normal in $D.$ Then by Lemma \ref{thm ZP}, there exist sequences $f_j\in \mathcal{F},$ $z_j \rightarrow z_0,$ $\rho_j \rightarrow 0,$ such that the sequence $\{g_j\}$ defined as  $g_j(\zeta )= \frac{f_j(z_j+\rho_j\zeta)}{\rho_j^{k}}$ converges  locally uniformly in $\mC^n$ to a non-constant entire function $g.$\\
 Let $w= (a_1,\ldots, a_n), ~w^{\prime}= (b_1,\ldots,b_n) \in\mC^n$ and define 
$$h_i(z_i):= g(b_1,\ldots,b_{i-1},z_i,a_{i+1},\ldots,a_n), ~i= 1,2,\ldots n.$$
By Hurwitz's theorem, it follows that $h_i(z)\neq 0 ~\mbox{for all}~ z\in\mC.$\\
By Weierstrass's theorem, we have  $$\frac{\partial^k}{\partial z_i^k}g_j(b_1,\ldots,b_{i-1},z_i,a_{i+1},\ldots,a_n) \to \frac{\partial^k}{\partial z_i^k}g(b_1,\ldots,b_{i-1},z_i,a_{i+1},\ldots,a_n).$$
\textbf{Claim:} $\frac{\partial^k}{\partial z_i^k}g(b_1,\ldots,b_{i-1},z_i,a_{i+1},\ldots,a_n) \not\equiv 1$

\smallskip

\noindent
Suppose $\frac{\partial^k}{\partial z_i^k}g(b_1,\ldots,b_{i-1},z_i,a_{i+1},\ldots,a_n) \equiv 1.$ Then $h_i(z)$ is a polynomial of degree $k$ and so $h_i(a) = 0$ for some $a\in\mC,$ a contradiction and this establishes the claim.

\smallskip

\noindent
Since $h_i(z)\neq 0 ~\mbox{for all}~ z\in\mC.$ So, $\frac{\partial^{k}}{\partial z_i^k}h_i(a)= 1 ~\mbox{for some}~ a\in\mC$ which implies 
$$\frac{\partial^k}{\partial z_i^k}g(b_1,\ldots,b_{i-1},a,a_{i+1},\ldots,a_n) = 1.$$
Since $\frac{\partial^k}{\partial z_i^k}g$ is non-constant. So, by Hurwitz's theorem there is a sequence $w_j = (b_1,\ldots,b_{i-1},a_j,a_{i+1},\ldots,a_n) \to w_0 = (b_1,\ldots,b_{i-1},a,a_{i+1},\ldots,a_n)$ such that for sufficiently large $j$
$$\frac{\partial^k}{\partial z_i^k}g_j(w_j) = 1$$ which implies
$$\sum\limits_{l=1}^{n} \frac{\partial^k }{\partial z_l^k}f_j(z_j+\rho_jw_j)e_l = 1$$

where $e_l=
\left\{\begin{array}{cc}
1, & ~\mbox{for}~ l=i\\
0, & ~\mbox{for}~ l\neq i.\end{array}\right.$\\
That is, $(D^kf_j(z_j+\rho_jw_j),v) = 1, ~\mbox{where}~ v = (0,\ldots,1,\ldots, 0),$ which is a contradiction and hence $\mathcal{F}$ is normal in $D.$  $\Box$

\medskip

The next Picard-type theorem is
\begin{theorem}\label{P1}
Let $f\in \mathcal{H}(\mC^n) ~\mbox{and}~ \ a ~\mbox{and}~ b$ be two distinct finite complex numbers. 
If $f\neq a$ and for each $\alpha\in f^{-1}(\{b\}) ,$ zero multiplicity of $f$ at $\alpha$ is at least two, then $f$ is constant. 
\end{theorem}
\textbf{Proof:}
 Let $w= (a_1,\ldots, a_n), ~w^{\prime}= (b_1,\ldots,b_n) \in\mC^n$ and define 
$$h_i(z_i)= f(b_1,\ldots,b_{i-1},z_i,a_{i+1},\ldots,a_n), ~i= 1,2,\ldots n.$$
By the given hypothesis $h_i$ omits $a.$ Suppose $h_i(\alpha_i) = b.$ Then 
$$f(b_1,\ldots,b_{i-1},\alpha_i,a_{i+1},\ldots,a_n) = b$$ and so $B_i = (b_1,\ldots,b_{i-1},\alpha_i,a_{i+1},\ldots,a_n)\in f^{-1}(\{b\}).$
Since for each $\alpha\in f^{-1}(\{b\}),$ zero multiplicity of $f$ at $\alpha$ is at least two, we have 
$$f(z) = \sum\limits_{m=2}^{\infty}P_m(z-B_i),$$ 
where $P_m$ is either identically zero or a homogeneous polynomial of degree $m.$\\
Now $$\frac{\partial}{\partial z_i}f(z) = \sum\limits_{m=1}^{\infty}q_m(z-B_i)$$
 where $q_m$ is either identically zero or a homogeneous polynomial of degree $m.$ \\
Thus $\frac{\partial}{\partial z_i}f(B_i) = 0$ which implies $\frac{d}{dz_i}h_i(\alpha_i) = 0.$ Thus $b-$points of $h_i$ have multiplicity at least two and so by Theorem \ref{thmN}, we conclude that each $h_i(z_i)$ is constant and hence $f$ is constant. $\Box$

\medskip

The corresponding normality criterion of Theorem \ref{P1} is 

\begin{theorem}\label{P2}
Let $D$ be a domain in $\mC^n$. Let $\mathcal{F}\subset \mathcal{H}(D), \ a ~\mbox{and}~ b$ be two distinct finite complex numbers.
If for each $f\in\mathcal{F}, ~f(z)\neq a$ and for each $\alpha\in f^{-1}(\{b\})\cap D,$ zero multiplicity of $f$ at $\alpha$ is at least two, then $\mathcal{F}$ is normal in $D.$  
\end{theorem}
\textbf{Proof:}
Suppose $\mathcal{F}$ is not normal in $D.$ Then by Lemma \ref{thm ZP}, there exist sequences $f_j\in \mathcal{F},$ $z_j \rightarrow z_0,$ $\rho_j \rightarrow 0,$ such that the sequence $\{g_j\}$ defined as  $g_j(\zeta )= f_j(z_j+\rho_j\zeta)$ converges  locally uniformly in $\mC^n$ to a non-constant entire function $g.$ Let $w= (a_1,\ldots, a_n), ~w^{\prime}= (b_1,\ldots,b_n) \in\mC^n$ and define 
$$h_i(z_i)= g(b_1,\ldots,b_{i-1},z_i,a_{i+1},\ldots,a_n), ~i= 1,2,\ldots n.$$
By Hurwitz's theorem, it follows that $h_i(z)\neq 0 ~\mbox{for all}~ z\in\mC.$\\
Suppose $h_i(\alpha_i) = b$ which implies $g(b_1,\ldots,b_{i-1},\alpha_i,a_{i+1},\ldots,a_n) = b.$ By Hurwitz's theorem there exist sequences $\beta_{j,i}= (b_1,\ldots,b_{i-1},\alpha_{j,i},a_{i+1},\ldots,a_n) \to \beta= (b_1,\ldots,b_{i-1},\alpha_i,a_{i+1},\ldots,a_n)$ such that $f_j(z_j+\rho_j\beta_{j,i}) = b$ which implies $z_j+\rho_j\beta_{j,i}\in f_j^{-1}(\{b\})\cap D,$ for sufficiently large $j.$ Since for $\alpha\in f_j^{-1}(\{b\})\cap D,$ zero multiplicity of $f_j$ at $\alpha$ is at least two, we have $f_j(z) = \sum\limits_{m=2}^{\infty}P_{j_m}(z-(z_j+\rho_j\beta_{j,i})),$ where $P_{j_m}$ is either identically zero or a homogeneous polynomial of degree $m.$\\
Thus $\frac{\partial}{\partial z_i}(f_j(z)) = \sum\limits_{m=1}^{\infty}q_{j_m}(z-(z_j+\rho_j\beta_{j,i})),$ where $q_{j_m}$ is either identically zero or a homogeneous polynomial of degree $m.$\\
Thus $$\frac{\partial}{\partial z_i}(f_j(z_j+\rho_j\beta_{j,i})) = 0$$ 
which implies $$\frac{\partial}{\partial z_i}g_j(b_1,\ldots,b_{i-1},\alpha_{j,i},a_{i+1},\ldots,a_n) = 0$$
Taking $j \to\infty, ~ \frac{\partial}{\partial z_i}g(b_1,\ldots,b_{i-1},\alpha_i,a_{i+1},\ldots,a_n) = 0$ which implies $\frac{d}{dz_i}h_i(\alpha_i)= 0.$ Thus $b-$ points of $h_i$ have multiplicity atleast two and so by  Theorem \ref{thmN} each  $h_i(z_i)$ is constant and hence $g$ is constant, a contradiction.           $\Box$

\medskip

The fundamental normality test for holomorphic functions of several complex variables is a simple consequence of Theorem \ref{P2}.

\section{Extension of Some Known Results from $\mC$ to $\mC^n$}
In $2012$, J. Grahl and S. Nevo\cite{Nevo} proved the following normality criterion:
\begin{theorem}\label{thm CN}
Let some $\epsilon> 0$ be given and set $\mathcal{F} := \{f\in\mathcal{H}(\Omega): f^{\sharp}(z)\geq \epsilon ~\mbox{for all}~ z\in \Omega\}.$ Then $\mathcal{F}$ is normal in $\Omega.$ 
\end{theorem}

We find a several complex variables analogue of Theorem \ref{thm CN}:

\medskip

\begin{theorem}\label{thmGN}
Let $D$ be a domain in $\mC^n$. Let $k\geq 0$ be an integer, $c>0, ~\alpha>1$ constant. Then 
$$\mathcal{F} = \left\{f\in\mathcal{H}(D): \sup\limits_{\|v\|=1}\frac{\lvert (D^kf(z),v)\rvert}{1+\lvert f(z)\rvert^{\alpha}}>c\right\}$$ is normal in $D.$
\end{theorem} 
\textbf{Proof:}
Suppose $\mathcal{F}$ is not normal in $D.$ Taking  $\beta>\frac{k}{\alpha-1},$  by Zalcman's lemma in $\mC^n,$  there exist sequences $f_j\in \mathcal{F},$ $z_j \rightarrow z_0,$ $\rho_j \rightarrow 0,$ such that the sequence $\{g_j\}$ defined as  
$$g_j(\zeta )=\rho_j^{\beta}f_j(z_j+\rho_j\zeta)$$ converges  locally uniformly in $\mC^n$ to a non-constant entire function $g.$ Let $\zeta_0\in\mC^n$ such that $g(\zeta_0)\neq 0.$ Clearly $f_j(z_j+\rho_j\zeta_0) \to \infty ~\mbox{as}~ j\to\infty.$ Also $$\sup\limits_{\|v\|=1}\lvert (D^kg_j(\zeta_0),v)\rvert =  \sup\limits_{\|v\|=1}\rho_j^{\beta+k}\lvert (D^kf_j(z_j+\rho_j\zeta_0),v)\rvert\to \sup\limits_{\|v\|=1}\lvert (D^kg(\zeta_0),v)\rvert.$$ 
Since
$$\sup\limits_{\|v\|=1}\lvert (D^kf_j(z_j+\rho_j\zeta_0),v)\rvert>c\lvert f_j(z_j+\rho_j\zeta_0)\rvert^\alpha,$$
we have 

\smallskip

\begin{eqnarray*}
 \sup\limits_{\|v\|=1}\lvert (D^kg_j(\zeta_0),v)\rvert&=&\sup\limits_{\|v\|=1}\rho_j^{\beta+k}\lvert (D^kf_j(z_j+\rho_j\zeta_0),v)\rvert\\
 &>&\rho_j^{\beta+k}c\lvert f_j(z_j+\rho_j\zeta_0)\rvert^\alpha\\
 &=& c(\rho_j^\beta\lvert f_j(z_j+\rho_j\zeta_0)\rvert)^\alpha. \rho_j^{\beta+k-\beta\alpha} \to\infty ~\mbox{as}~ j\to \infty.
\end{eqnarray*}
which is a contradiction to the fact that $g$ is  holomorphic at $\zeta_0.$ Therefore, $\mathcal{F}$ has to be normal. $\Box$

\begin{remark}
For $k = 1 ~\mbox{and}~ \alpha = 2,$ we obtain several complex variable version of Theorem \ref{thm CN}.\\
\end{remark}

Further, we obtain an improvement of a Marty's Theorem in $\mC^n$:
\begin{theorem}\label{Marty}
Let $D$ be a domain in $\mC^n,$  $\mathcal{F}\subset \mathcal{H}(D)$ and $\alpha$ be a positive real number. If 
$$F_{\alpha} = \left\{\sup\limits_{\|v\|=1}\frac{\lvert (Df(z),v)\rvert}{1+\lvert f(z)\rvert^{\alpha}}:f\in\mathcal{F}\right\}$$ is locally uniformly bounded in $D,$ then $\mathcal{F}$  is normal in $D.$ 
\end{theorem}
\textbf{Proof:} 
Suppose $\mathcal{F}$ is not normal in $D.$ Then by Lemma \ref{thm ZP}, there exist sequences $f_j\in \mathcal{F},$ $z_j \rightarrow z_0,$ $\rho_j \rightarrow 0,$ such that the sequence $\{g_j\}$ defined as  $g_j(\zeta )= \rho_j^{-1/2}f_j(z_j+\rho_j\zeta)$ converges  locally uniformly in $\mC^n$ to a non-constant entire function $g.$\\
By Weierstrass' Theorem, $$\frac{\partial}{\partial z_i}g_j(\zeta)\to \frac{\partial}{\partial z_i}g(\zeta), ~i = 1,2,\ldots,n.$$\\
Now 
\begin{eqnarray*}
\sup\limits_{\|v\|=1}\sum_{l=1}^{n}\left\lvert \frac{\partial}{\partial z_l}g_j(\zeta).v_l\right\rvert &=& \rho_j^{1/2}\sup\limits_{\|v\|=1}\sum_{l=1}^{n}\left\lvert \frac{\partial}{\partial z_l}f_j(z_j+\rho_j\zeta).v_l\right\rvert\\
                           &=& \rho_j^{1/2} \sup\limits_{\|v\|=1}\lvert (Df_j(z_j+\rho_j\zeta),v)\rvert\\
													 &\leq& \rho_j^{1/2} M[1+\lvert f_j(z_j+\rho_j\zeta)\rvert^{\alpha}]\\
													&=& \rho_j^{1/2} M[1+\lvert \rho_j^{1/2}g_j(\zeta)\rvert^{\alpha}] \to 0 ~\mbox{as}~ j\to \infty
\end{eqnarray*}
Thus 
$$\sup\limits_{\|v\|=1}\sum_{l=1}^{n}\left\lvert \frac{\partial}{\partial z_l}g_j(\zeta).v_l\right\rvert = 0.$$
That is,  $$\frac{\partial g}{\partial z_i}(\zeta) = 0, ~i= 1,\ldots,n,$$
which  implies $g$ is a constant, a contradiction.         $\Box$

\medskip

\bibliography{sn-bibliography}


\end{document}